# Restrictions on collapsing with a lower sectional curvature bound

Vitali Kapovitch

October 29, 2018

## Abstract

We obtain new topological information about the local structure of collapsing under a lower sectional curvature bound. As an application we prove a new sphere theorem and obtain a partial result towards the conjecture that not every Alexandrov space can be obtained as a limit of a sequence of Riemannian manifolds with sectional curvature bounded from below.

## 1 Introduction

The study of Alexandrov spaces with curvature bounded from below has been largely motivated by the fact that they naturally appear as the boundary points of the class of Riemannian manifolds with sectional curvature bounded from below.

If $X$ is a limit of a sequence of Riemannian manifolds $M_i^n$ satisfying $\sec \geq k$, then it is not hard to see [BGP92] that the Hausdorff dimension of $X$ can not be greater than $n$.

If it is equal to $n$ we say that the sequence $M_i^n$ converges without collapse and if it is less than $n$ we say that this sequence collapses.

The first case is understood fairly well at least topologically due to the stability theorem by Perelman [Per91], which says that for sufficiently large indices, Hausdorff approximations $M_i \to X$ are close to homeomorphisms if $X$ is compact.

Unlike the situation in the noncollapsing case, very little is known about the structure of the limit and its relationship to the elements of the sequence when collapse does occur. The main structural result here is Yamaguchi's Fibration Theorem [Yam91] which asserts that if the limit is a Riemannian manifold, then the Hausdorff approximations into the limit from the elements of the sequence can be chosen to be smooth fibrations. However, practically nothing is known about the structure of collapse when the limit space is singular.

In the present paper we obtain at least some partial understanding of the local topology of collapsing when the limit is an *arbitrary* Alexandrov space. We prove:

**Theorem 1.1.** *Let $k \in \mathbb{R}$ and suppose $M_m^n$ is a sequence of Riemannian manifolds with $\sec \geq k$, Gromov-Hausdorff converging to an Alexandrov space $X$.*

*Then for any $x_0 \in X$, there exists an $r_0 = r_0(x_0) > 0$ such that for any sequence of points $x_m \in M_m$ converging to $x_0$, for any sufficiently large $m$, the closed ball $\bar{B}(x_m, r_0)$ is a manifold with boundary simply homotopy equivalent to a finite CW complex of dimension $\leq n - \dim X$.*

Let $\bar{\mathcal{M}}_k^n$ be the closure in the pointed Gromov-Hausdorff topology of the class of $n$-dimensional Riemannian manifolds with $\sec \geq k$.

The following natural question remains unanswered:

**Question 1.2.** *Is it true that for any finite-dimensional Alexandrov space $X$, there exist $k$ and $n$ such that $X \in \bar{\mathcal{M}}_k^n$?*

The collapsing phenomenon occurs naturally when one considers a pointed sequence formed by rescaling of a nonnegatively curved open manifold by positive constants approaching 0. The limit in this case is a Euclidean cone over the ideal boundary of $M$, which we will denote by $M(\infty)$.

From this description it is easy to conclude that $M(\infty)$ is an Alexandrov space with curvature bounded below by 1 [BGP92]. It was shown in [GK95] that if the ideal boundary is a Riemannian manifold, then its topology is severely restricted:

**Theorem 1.3.** *[GK95] Let $M^n$ be a complete open manifold with $K_{sec} \geq 0$. If $M(\infty)$ is an $m$-dimensional Riemannian manifold, then it has a finite covering $N \to M(\infty)$ satisfying one of the following conditions:*

(i)   *$N$ is homotopy equivalent to $S^m$, or*

(ii)  *$N$ is homotopy equivalent to $\mathbb{CP}^{m/2}$, or*

(iii) *$N$ has rational cohomology ring of $\mathbb{HP}^{m/4}$*

It was also shown in [PWZ95] that if the spherical suspension over a Riemannian manifold $N$ belongs to $\bar{\mathcal{M}}_{1/4}^n$, then $N$ must satisfy one of the conditions $(i) - (iii)$ above.

Based on this and Theorem 1.3, it was conjectured in [GK95] that the answer to Question 1.2 is negative and, more specifically, that if $M$ is a positively curved manifold that does not satisfy the conclusion of Theorem 1.3, then the spherical suspension $\mathbb{S}M$ does not belong to $\bar{\mathcal{M}}_k^n$ for any $k, n$.

We use Theorem 1.1 to obtain a partial result towards verifying this conjecture.

Given an Alexandrov space $X$, we will say that the minimal collapsing codimension of $X$ is equal to $s$ if $s$ is the smallest integer such that $X \in \bar{\mathcal{M}}_k^{\dim X + s}$ for some $k \in \mathbb{R}$. If no such $s$ exists we will say that the minimal collapsing codimension of $X$ is equal to $\infty$.

**Example 1.4.** The minimal collapsing codimension of the spherical suspension over $CP^n$ is equal to 1 if $n > 1$. Indeed, by Perelman's stability theorem, it has to be positive. On the other hand, as was observed by Yamaguchi [Yam91], by rescaling the fibers of the standard $S^1$ action on $S(S^{2n+1}) = S^{2n+2}$, one can find a sequence of metrics with $\sec \geq 0$ on $S^{2n+2}$ collapsing to $\mathbb{S}CP^n$.

It is well-known (cf. [Ber61]) that the 24-dimensional Caley flag $\Sigma^{24} = F_4/Spin(8)$ admits a homogeneous metric of $sec \geq 1$.

We prove

**Theorem 1.5.** *Let $X$ be an Alexandrov space such there exists a point $x_0 \in X$ such that the space of directions to $X$ at $x_0$ is a Riemannian manifold.*

(a) *If $\Sigma_{x_0} X$ is diffeomorphic to $\Sigma^{24}$, then the minimal collapsing codimension of $X$ is $\geq 15$;*

(b) *If $\Sigma_{x_0} X$ is diffeomorphic to the Caley plane $CaP^2$, then the minimal collapsing codimension of $X$ is $\geq 8$;*

(c) *If $\Sigma_{x_0} X$ is diffeomorphic to $\mathbb{H}P^n$, then the minimal collapsing codimension of $X$ is $\geq 3$.*

**Remark 1.6.** If we equip $\Sigma^{24}$ with the Berger metric of $sec \geq 1$ then the natural suspension metric of $curv \geq 1$ on $X = S\Sigma$ satisfies assumption (a) of Theorem 1.5. Similarly, spherical suspensions $\mathbb{S}Cap^2, \mathbb{S}\mathbb{H}P^n$ over the symmetric spaces $CaP^2$ and $\mathbb{H}P^n$ of $sec \geq 1$ satisfy assumptions of (b) and (c) respectively.

Observe that unlike the Caley plane and the Caley flag $F_4/Spin(8)$, the quaternionic projective space does satisfy the conclusion of Theorem 1.3. Moreover, the same construction as in Example 1.4 shows that there exists a sequence of metrics with $sec \geq 0$ on $S^{4n+4}$ Gromov-Hausdorff converging to $\mathbb{S}\mathbb{H}P^n$.

Hence, the bound provided by theorem 1.5 in this case is *sharp* and the minimal collapsing codimension of $\mathbb{S}\mathbb{H}P^n$ is equal to 3. Furthermore, the following sphere theorem shows that in some sense the above example of collapsing to $\mathbb{S}\mathbb{H}P^n$ is the only one possible if the codimension of collapse is equal to 3:

**Theorem 1.7.** *Let $k \in \mathbb{R}$ and $n > 1$ be an integer. There exists an $\epsilon = \epsilon(k,n) > 0$ such that if $M^{4n+4}$ is a complete Riemannian manifold satisfying $\sec(M^{4n+4}) \geq k$ and $d_{G-H}(M, \mathbb{S}\mathbb{H}P^n) \leq \epsilon$, then $M^{4n+4}$ is homeomorphic to $S^{4n+4}$.*

Let us briefly describe the strategy of the proofs.

To prove Theorem 1.1, we use a special kind of averaging procedure for distance functions due to Perelman [Per93] (cf. [PP93], [Kap99]) to construct a strictly convex function $f$ near $x_0 \in X$ with a minimum at $x_0$ and then lift it to the elements of the sequence. Standard critical point theory for distance functions implies that the sublevel sets of the lifts are homeomorphic to closed balls $\bar{B}(x_m, r_0)$ if $r_0 = r_0(x_0)$ is sufficiently small and $d_{G-H}(M_m, X) \ll r_0$.

Our crucial observation is that the function $f$ can be chosen in such a way that the lifts $f_m$ are partially convex in the sense of Wu. We will give a careful definition of partial convexity in Section 3, but informally speaking, we will show that the sum of any $s+1$ eigenvalues of the "hessian" of $f_m$ at any point is positive.

Using the approximation result of Wu [Wu87], we can assume that $f_m$ is smooth and hence its hessian at any point has at most $s$ nonpositive eigenvalues. By further approximating $f_m$ by a Morse function, we conclude that the sublevel sets of $f_m$ are homeomorphic to sublevel sets of a Morse function with indices of critical points at most $s$ which immediately yields the statement of Theorem 1.1.

To prove Theorem 1.5, we observe that, by a standard rescaling argument, we can assume that $X$ is isometric to the Euclidean cone over $\Sigma$ with $x_0$ equal to the vertex. Since $X$ is smooth away from $x_0$, we can use Yamaguchi's fibration theorem to conclude that metric spheres at $x_m$ fiber over $\Sigma$ with closed manifolds as the fibers. Theorem 1.1 imposes certain obvious restriction on the cohomology of the metric spheres centered at $x_m$. We then use a Serre spectral sequence argument to show that a total space of a bundle over $\Sigma_{x_0} X$ can not satisfy these restrictions if the dimension of the fibers is too small.

For the proof of Theorem 1.7, we observe that by a standard critical point theory argument, $M_m$ is homeomorphic to a union of two metric balls of fixed radius satisfying the conclusion of Theorem 1.1, glued along a homeomorphism of the boundaries. Arguing as in the proof of Theorem 1.5, we see that the boundaries of these balls fiber over $\mathbb{HP}^n$. We then show that in order to satisfy the homological restriction implied by Theorem 1.1, the metric spheres in question must be homeomorphic to $S^{4n+3}$. Again using Theorem 1.1, we then conclude that the corresponding metric balls are contractible. By Poincare conjecture, this immediately implies that $M_m$ is homeomorphic to a sphere.

**Acknowledgements.**

The author would like to thank Igor Belegradek, Burkhard Wilking and Wolfgang Ziller for helpful conversations and suggestions related to this paper.

## 2   Notations and conventions

Throughout this paper all homology and cohomology groups have $\mathbb{Z}$ coefficients unless otherwise indicated.

For an Alexandrov space $X$ we will denote by $CX$ the Euclidean cone over $X$ and by $\mathbb{S}X$ the spherical suspension over $X$. We will denote the spherical join of $X$ and $Y$ by $X * Y$. For a point $x$ in an Alexandrov space $X$ we will denote the space of directions of $X$ at $x$ by $\Sigma_x X$. The reader is referred to [BGP92] or [BBI01] for the definition of a space of directions and other basic notions of Alexandrov geometry.

Let $p, q \in X$ be two points in a finite dimensional Alexandrov space $X$. We will use the following notation:

$p(q)' = \{\xi \in \Sigma_q X|$ there exists a shortest geodesic $\gamma$ from $q$ to $p$ such that $\gamma'(0) = \xi\}$. Observe that $p(q)'$ is always closed. With this notation the first variation formula takes the following form

$$d(\cdot, p)'(\xi) = -\cos \angle \xi p(q)'$$

for any $\xi \in \Sigma_q X$.

Let $f$ be a Lipschitz function on a Riemannian manifold $M$. Let $V$ be a $C^\infty$ vector field on $M$ and let $U$ be any subset of $M$. We will say that $V$ is *gradient-like* for $f$ on $U$ if the directional derivative $f'(V)$ exists everywhere in $U$ and moreover there exists a constant $c > 0$ such that $f'(V(x)) \geq c$ for any $x \in U$.

## 3 Partially concave functions

The notion of partially convex functions was introduced by H. Wu in [Wu87]. In this paper we will work with the dual notion of partially concave functions. For convenience of the reader we will reproduce the relevant definitions.

Let $f\colon M \to \mathbb{R}$ be a continuous function on a Riemannian manifold $M$. Let $\gamma\colon (-a, a) \to M$ be a geodesic such that $\gamma(0) = x \in M$ and $\gamma'(0) = X \in T_x M$.

Define

$$Cf(x; X) = \limsup_{r \to 0} \frac{1}{r^2}\{f(\gamma(r)) + f(\gamma(-r)) - 2f(\gamma(0))\} \tag{3.1}$$

We say that a set of $s$ vectors $\{X_1, \ldots, X_s\}$ in an inner product space $V$ is $\epsilon$-orthonormal if $|\langle X_i, X_j\rangle - \delta_{i,j}| < \epsilon$ for all $i, j$.

Let $M^n$ be a Riemannian manifold and let $s \leq n$ be a positive integer.

**Definition 3.1.** *We say that a function $f\colon M \to \mathbb{R}$ belongs to the class $\mathcal{T}(s)$ if $f$ is locally Lipschitz and for each $x_0 \in M$ there exists a neighborhood $W$ of $x_0$ and constants $\epsilon, \eta > 0$ such that*

$$\sum_{i=1}^{s} Cf(x, X_i) \leq -\eta$$

*for any $x \in W$ and $X_1, \ldots, X_s$ - an $\epsilon$-orthonormal set in $T_x M$.*

Note that $\mathcal{T}(1)$ is equal to the set of all strictly concave functions on $M$ and as it was shown in [Wu87], $\mathcal{T}(n)$ is the set of all locally Lipschitz strictly superharmonic functions on $M$.

**Remark 3.2.** It is immediate to check that a positive linear combination and the minimum of a finite number of functions from $\mathcal{T}(s)$ again belongs to $\mathcal{T}(s)$.

We will make use of the following approximation result proved in [Wu87]:

**Theorem 3.3.** *Let $f \in \mathcal{T}(s)$ where $1 \leq s \leq \dim M$ and let $\epsilon\colon M \to \mathbb{R}$ be a positive continuous function.*

*Then there exists a $C^\infty$ function $F \in \mathcal{T}(s)$ such that $|F - f| < \epsilon$.*

## 4 Morse theory for partially convex functions

The proof of the following well-known Lemma is an elementary exercise in basic algebraic topology.

**Lemma 4.1.** *Let $W^{n+1}$ be a compact oriented manifold which is a thickening of a $s$-dimensional CW complex. Then*

$$H^i(\partial W, A) = H_i(\partial W, A) = 0$$

*for any ring of coefficients $A$ and any $i$ satisfying $s < i < n - s$.*

*Sketch of the proof.* Suppose $A = \mathbb{Z}$. Since $W^{n+1}$ is homotopy equivalent to an $s$-dimensional complex, $H^i(W) = H_i(W) = 0$ for any $i > s$. By Poincare duality, this implies that $H^i(W, \partial W) = H_i(W, \partial W) = 0$ for $i < n + 1 - s$. Now the claim of the lemma immediately follows from the long exact homology and cohomology sequences of the pair $(W, \partial W)$. The case of general $A$ follows from the case $A = \mathbb{Z}$ by the universal coefficients formula. □

The fundamental theorem of Morse theory implies that if $f \in \mathcal{T}(s+1)$ is a Morse function on $M^n$ with compact superlevel sets, then for any regular value $c$ of $f$, the superlevel set $\{f \geq c\}$ can be obtained from $\emptyset$ by attaching a finite number of handles of index at most $s$ (or equivalently, $\{f \geq c\}$ can be obtained from $\{f = c\}$ by attaching a finite number of handles of index at least $n - s$ ).

The next Lemma shows that it is also true for arbitrary functions from $\mathcal{T}(s+1)$ once the notion of a regular value is properly understood.

**Lemma 4.2.** *Let $h \colon M^{n+1} \to \mathbb{R}$ be a function from $\mathcal{T}(s+1)$ with compact superlevel sets on an orientable manifold $M$. Let $[c_1, c_2] \subset \mathrm{Im}(h)$. Suppose that the following two conditions are satisfied*

(i) *$h$ has directional derivatives everywhere in $h^{-1}([c_1, c_2])$ and moreover there exists an $L > 0$ such that the derivative $h'_x$ is $L$-Lipschitz on $T_x M$ for any $x \in h^{-1}([c_1, c_2])$*

(ii) *there exists a gradient-like smooth vector field $X$ for $h$ on $h^{-1}([c_1, c_2])$.*

*Then there exists a Morse function $\hat{h}$ uniformly close to $h$ on $h^{-1}([c_1, \infty))$ such that*

1. *any $c \in [c_1, c_2]$ is a regular value of $\hat{h}$*
2. *$\hat{h} \in \mathcal{T}(s+1)$ on $\{h \geq c_1\}$*
3. *$\{h \geq c\}$ is homeomorphic to $\{\hat{h} \geq c\}$ for any $c \in [c_1, c_2]$.*

**Remark 4.3.** Lemma 4.1 immediately implies that under the assumptions of Lemma 4.2,

$$H^i(\{\hat{h} = c\}, A) = H_i(\{\hat{h} = c\}, A) = 0$$

for any ring of coefficients $A$ and any $i$ satisfying $s < i < n - s$.

*Proof of Lemma 4.2.* The proof is essentially an application of the Smoothing Theorem of Wu mentioned in Section 3. Unfortunately, the result we want follows from the proof rather than the statement of that theorem. Therefore we will briefly outline the construction involved in its proof.

Let $\kappa \colon R \to [0,1]$ be a $C^\infty$ function with support in $[-1,1]$ such that $\kappa = const$ near 0 and
$$\int_{\mathbb{R}^n} \kappa(|v|) dv = 1$$

Define $h_\rho \colon M \to \mathbb{R}$ by the formula

$$h_\rho(x) = \frac{1}{\rho^n} \int_{T_xM} f(exp_x v) \kappa(\rho v) d\mu_v$$

where $d\mu_v$ stands for the Lebesgue measure on $T_xM$.

Then according to the proof of [Wu87, Lemma 2], $h_\rho$ is smooth and belongs to $\mathcal{T}(s+1)$ on $h^{-1}([c_1, +\infty))$ if $\rho$ is sufficiently small. Fix a $c \in (c_1, c_2)$ We are going to show that if $\rho$ is sufficiently small then the superlevel sets of $h_\rho$ and $h$ are homeomorphic.

First of all let us show that $c$ is a regular value of $h_\rho$.

According to [GS77, Proposition 2.1], the differential of $h_\rho$ can be computed as follows:

Let $u \in T_xM$.

Construct a vector field $U$ on $B_\rho(x)$ as follows.

Let $\gamma$ be the unique geodesic with $\gamma'(0) = u$ and define for each $y \in B_\rho(x)$ a smooth curve $\gamma_y$ by the formula

$$\gamma_y(t) = \exp_{\gamma(t)}(P_{\gamma(t)}(\exp^{-1}_{\gamma(0)}(y)))$$

Observe that $U$ is well defined and smooth if $\rho < \mathrm{injrad} M$.

Then $dh_\rho(u)$ is given by the following formula

$$dh_\rho(u) = \frac{1}{\rho^n} \int_{T_xM} h'_{exp(v)}(U) \kappa(\rho v) d\mu_v \tag{4.2}$$

Let $u = X(x)$. By construction of $U$, we see that $U$ is close to $X$ on $B_\rho(x)$ if $\rho$ is sufficiently small, which by the Lipschitz condition (i) on $h'$ implies that $dh_\rho(X(x))$ is uniformly close in $h^{-1}([c_1, +\infty))$ to

$$\frac{1}{\rho^n} \int_{T_xM} h'_{exp(v))}(X(\exp(v)) \kappa(\rho v) d\mu_v$$

Now condition (ii) on $h$ implies that $X$ is a gradient-like vector field for $h_\rho$ on $h^{-1}([c_1, c_2])$ for all sufficiently small $\rho$.

Since $X$ is gradient-like for both $h$ and $h_\rho$, a standard argument using the flow of $X$ implies that $\{h \geq c\}$ and $\{h_\rho \geq c\}$ are homeomorphic for all sufficiently small $\rho$. Since Morse functions are dense among $C^\infty$ functions in the $C^\infty$ topology, we can assume that $h_\rho$ is Morse which concludes the proof of Lemma 4.2. □

# 5 Concavity of distance functions on Alexandrov spaces

In [Per93] Perelman introduced the following definition

**Definition 5.1.** *A function $f\colon U \to \mathbb{R}$ defined on a domain $U$ in an Alexandrov space $X$ is called $\lambda$-concave if for any unit speed shortest geodesic $\gamma \subset U$ the function $t \mapsto f(\gamma(t)) + \lambda t^2$ is concave.*

Observe that a Lipschitz function on $X$ is $\lambda$-concave iff $Cf(x;v) \leq -\lambda$ for any $x \in X, v \in \Sigma_x X$.

It is trivial to check that a positive linear combination or the infimum of a family of $\lambda$-concave functions is $\lambda$-concave. Also, a pointwise limit of a sequence of $\lambda$-concave functions is again $\lambda$-concave.

Toponogov triangle comparison implies that distance functions in a space of curvature $\geq k$ are more concave than distance functions in the model space of constant curvature $k$ and therefore it is easy to see that the following property holds [Per93]:

Let $p, q \in X$ be two points in an Alexandrov space $X$ of $curv \geq k$. Let $d = d(p,q)$ and $\epsilon < d/2$. Then $f(\cdot) = d(\cdot, q)$ is $\lambda$-concave in $B(p, \epsilon)$ where $\lambda$ depends *only* on $d$ and the lower curvature bound $k$.

**Remark 5.2.** The class of examples of $\lambda$-concave functions given by the distance functions can be enlarged using the following simple but important observation from [Per93]: If $f$ is $\lambda$-concave with $\lambda < 0$ and $\phi\colon \mathbb{R} \to \mathbb{R}_+$ is a concave $C^2$ function satisfying $0 \leq \phi' \leq 1$ then $\phi(f)$ is again $\lambda$-concave. Indeed, it is clearly enough to consider $f\colon \mathbb{R} \to \mathbb{R}$. If $f$ is $C^2$ then $\lambda$-concavity of $f$ is equivalent to the inequality $f'' \leq -\lambda$. Computing the second derivative of $\phi(f)$ we observe:

$$\phi(f)'' = \phi''(f)(f')^2 + \phi'(f)f'' \leq \phi'(f)f'' \leq \phi'(f)(-\lambda) \leq -\lambda$$

The general case immediately follows from this one since any $\lambda$-concave function on $\mathbb{R}$ can be approximated by $C^\infty$ $\lambda$-concave functions.

# 6 Local topology of collapsed spaces

In this section we prove Theorem 1.1 stated in the introduction. In fact we are going to prove the following more general statement.

**Theorem 6.1.** *Suppose $M_m^n$ is a sequence of Riemannian manifolds with $\sec \geq k$, Gromov-Hausdorff converging to an Alexandrov space $X$. Let $x_0 \in X$. Then there exists an $r_0 = r_0(x_0) > 0$ such that for any sequence of points $x_m \in M_m$ converging to $x_0$, for any sufficiently large $m$, the closed ball $\bar{B}(x_m, r_0)$ is homeomorphic to a smooth compact manifold with boundary $W_m^n$ which can be obtained from $\emptyset$ by attaching a finite number of handles of index at most $n - \dim X$.*

*Proof.* First of all we will reduce the situation to the case when $X = T_{x_0}$. Let $\epsilon_m = d_{G-H}(X, M_m)$. An easy rescaling argument [Kap99, Lemma 3.1] shows that that

$$(\frac{1}{\sqrt{\epsilon_m}} M_m, x_m) \to (T_{x_0} X, o)$$

Fix a sufficiently small $R$ such that $\frac{1}{R} B(x_0, R)$ is Hausdorff close to $B(o, 1)$ in $T_{x_0} X$.

A standard comparison argument (cf. Lemma 6.2 below) shows that for any fixed $r$, the function $d(\cdot, x_m)$ has no critical points in $\bar{B}(x_m, R/2) \backslash B(x_m, \sqrt{\epsilon_m} r)$ if $m$ is sufficiently large.

Suppose we can prove that there exists an $r_0 > 0$ such that $r_0$-balls around $x_m$ in $\frac{1}{\epsilon_m} M_m$ satisfy the conclusion of Theorem 6.1. Notice that $B_{\frac{1}{\sqrt{\epsilon_m}} M_m}(x_m, r_0) = B_{M_m}(x_m, \sqrt{\epsilon_m} r_0)$ which by above is homeomorphic to $B(x_m, R)$ and hence, the conclusion of Theorem 6.1 holds for $B(x_m, R)$ as well.

Let us therefore from now on assume that $X = T_{x_0} X$ and $(M_m, x_m) \to (T_{x_0} X, o)$ to begin with.

Let $\mu_m = d_{G-H}(X, M_m)$ and $h_m \colon X \to M_m$ be a $\mu_m$-Hausdorff approximation. Let $x_m = h_m(x_0)$. Our goal is to prove that there exists an $r > 0$ such that $\bar{B}(x_m, r)$ is a thickening of a $CW$-complex of dimension $\leq s = n - \dim X$ for any sufficiently large $m$.

The proof of the following Lemma is an elementary exercise in Toponogov angle comparison.

**Lemma 6.2.** *Let $0 < a < b$ be fixed constants. Then for any sufficiently large $m$ there exists a $C^\infty$ unit vector field $V_m$ on the annulus $\bar{B}(x_m, b) \backslash B(x_m, a)$ satisfying $d(\cdot, x_m)'(V_m) \geq 1 - O(\mu_m)$.* □

To produce $f$ and $f_m$ we use the same construction as in [Kap99, Theorem 1.3].

Let $\{q_\alpha\}_{\alpha=1}^N$ be a maximal $\pi/16$-separated net in $\Sigma = \Sigma_{x_0} X$. A standard argument shows that the balls $B(q_\alpha, \pi/8)$ cover $\Sigma$.

Fix a small $\delta > 0$. Throughout the rest of the proof we will denote by $c_i$ or $c$ various positive constants depending on $n, \Sigma$ but not on $\delta$.

For each $\alpha$, choose a collection $\{q^{\alpha\beta}\}_{\beta=1}^{N_\alpha}$ to be a maximal $\delta$-separated net in $B(q^\alpha, \pi/16)$ where the ball is taken inside $\Sigma$. Let $d = \dim X$.

A standard volume comparison argument shows that $N_\alpha$ satisfies

$$c_1/\delta^{d-1} \geq N_\alpha \geq c_2/\delta^{d-1} \tag{6.3}$$

for any $\alpha$.

Let $\phi_\delta : \mathbb{R} \to \mathbb{R}$ be the continuous function uniquely determined by the following properties:

(1) $\phi_\delta(0) = 0$

(2) $\phi'_\delta(t) = 1$ for $t \leq 1 - \delta^3$

(3) $\phi'_\delta(t) = 1/2$ for $t \geq 1 + \delta^3$

(4) $\phi''_\delta(t) = -1/(4\delta^3)$ for $1 - \delta^3 < t < 1 + \delta^3$

For every $\alpha$ define $f_\delta^\alpha$ by the following formula:

$$f_\delta^\alpha(x) = \frac{1}{N_\alpha} \sum_{\alpha=1}^{N_\alpha} \phi_\delta(d(x, q^{\alpha\beta}))$$

Put

$$f_\delta = \min_\alpha f_\delta^\alpha$$

Then according to Lemma 3.6 from [Per93], (cf. [PP93, Lemma 4.3]), the function $f_\delta^\alpha$ (and hence $f_\delta$) is strictly concave in $B(o, \delta^3/2)$ for all sufficiently small $\delta$.

Observe that $f_\delta(o) = \phi_\delta(1)$ since $d(o, q^{\alpha\beta}) = 1$ for any $\alpha\beta$. Moreover, we claim that $o$ is a point of a strict local maximum of $f_\delta$. Indeed, let $x \in X$ be a point sufficiently close to $o$. Without too much abuse of notations we can write $x = t\xi$ for some $\xi \in \Sigma$.

Since $\cup_\alpha B(q_\alpha, \pi/8) = \Sigma$, there exists $\alpha_0$ such that $\angle \xi q^{\alpha_0} \leq \pi/8$. Therefore, $\angle \xi q^{\alpha_0\beta} \leq 3\pi/16$ for any $\beta$ and hence, by the Toponogov angle comparison, $d(t\xi, q^{\alpha_0\beta}) \leq 1 - 2\cos(3\pi/16)t$ for any $\beta$ and any $t \leq c_0$.

By monotonicity of $\phi_\delta$, this implies that there exist a universal $\eta > 0$ such that

$$f_\delta^{\alpha_0}(t\xi) \leq f_\delta(o) - t\eta$$

if $t \leq c_0$ and since $f_\delta(t\xi) \leq f_\delta^{\alpha_0}(t\xi)$ we finally obtain that

$$f_\delta(t\xi) \leq f_\delta(o) - t\eta$$

for all $\xi \in \Sigma$ and $t \leq c_0$.

Since $\eta$ is fixed we can from now on assume that $\delta \ll \eta$.

By continuity of $f_\delta$, there exists $\nu = \nu(\delta) \ll \delta^3$ such that

$$\inf_{x \in S(o, \nu\delta^3)} f_\delta(x) > \sup_{x \in S(o, \delta^3/4)} f_\delta(x) \tag{6.4}$$

This implies that there exists a positive constant $a$, such that the level set $\{f_\delta = a\}$ is entirely contained in the open annulus $\{x \in X | \nu\delta^3 < d(x, o) < \delta^3/4\}$ i.e

$$\{f_\delta = a\} \subset B(o, \delta^3/4) \backslash \bar{B}(o, \nu\delta^3) \tag{6.5}$$

Let us lift $f_\delta$ to the elements of the sequence $(M_m, x_m)$ in a natural way. That is, put $q_m^\alpha = h_m(q^\alpha)$, $q_m^{\alpha\beta} = h_m(q^{\alpha\beta})$ and put

$$f_\delta^{\alpha m}(y) = \frac{1}{N_\alpha} \sum_{\beta=1}^{N_\alpha} \phi_\delta(d(y, q_m^{\alpha\beta}))$$

$$f_\delta^m = \min_\alpha f_\delta^{\alpha m}$$

Let us examine this function more carefully. First of all, notice that (6.5) implies that $f_\delta^m$ has compact superlevel sets and

$$\{f_\delta^m = a\} \subset B(x_m, \delta^3/4) \backslash \bar{B}(x_m, \nu\delta^3) \tag{6.6}$$

for all sufficiently large $m$.

By the definition, $f_\delta^m$ is Lipschitz and moreover the first variation formula implies that it has directional derivatives everywhere in $B(x_m, 1/2)$.

By lemma 6.2, the distance function $d(\cdot, x_m)$ has no critical points in the annulus $\bar{B}(x_m, \delta^3/4) \backslash B(x_m, \nu\delta^3)$ if $m$ is sufficiently large. The following lemma shows that the same remains true for $f_\delta^m$.

Let $V_m$ be the almost radial smooth vector field on $\bar{B}(x_m, \delta^3/4) \backslash B(x_m, \nu\delta^3)$ whose existence is guaranteed by Lemma 6.2.

**Lemma 6.3.** *For all sufficiently large $m$, $(f_\delta^m)'(-V_m(y)) > c$ and $(f_\delta^m)'(V_m(y)) < -c$ for any $y \in \bar{B}(x_m, \delta^3/4) \backslash B(x_m, \nu\delta^3)$.*

*Proof.* We can assume that $m$ is big enough so that $\mu_m \ll \delta^3$. Let $y_m \in \bar{B}(x_m, \delta^3/4)$. By the chain rule we see that for any $v \in \Sigma_{y_m} M_m$ and any $\alpha$.

$$(f_\delta^{\alpha m})'(v) = \frac{1}{N} \sum_{\alpha=1}^{N_\alpha} \phi_\delta'(d(y, q_m^{\alpha\beta})) d(y, q_m^{\alpha\beta})'(v) \tag{6.7}$$

Since $h_m$ is a $\mu_m$- Hausdorff approximation, there is $y \in X$ such that $d(y_m, h_m(y)) < \mu_m$. As before we will write $y$ as $y = t\xi$. Let $\xi_m = h_m(\xi)$.

Let $\alpha_0$ be such that

$$f_\delta^m(y_m) = f_\delta^{\alpha_0 m}(y_m) \tag{6.8}$$

It is easy to check that $q^{\alpha_0}$ must satisfy

$$\angle \xi q^{\alpha_0} \leq \pi/4 \tag{6.9}$$

and therefore

$$\angle \xi q^{\alpha_0 \beta} \leq 5\pi/16 \tag{6.10}$$

for any $\beta$.

Since $\mu_m \ll \delta$, we see that $|\tilde\angle x_m y_m q_m^{\alpha_0 \beta} - \tilde\angle oyq^{\alpha_0\beta}| \leq O(\mu_m)$. On the other hand $\tilde\angle oyq^{\alpha_0\beta} = \angle oyq^{\alpha_0\beta} = \pi - \angle yoq^{\alpha_0\beta} - \angle oq^{\alpha_0\beta}y \geq \pi - 5\pi/16 + \eta - O(\delta) = 9\pi/16 - O(\delta)$ by (6.10). Since $\delta \ll \eta$ we can conclude that $\tilde\angle oyq^\alpha \geq 17\pi/32$ and therefore $\tilde\angle x_m y_m q_m^{\alpha_0\beta} \geq 33\pi/64$ if $m$ is sufficiently large. By the Toponogov angle comparison, this means that $\angle vu \geq 33\pi/64$ for any $v \in x_m(y_m)', u \in q_m^{\alpha_0\beta}(y_m)'$

By the first variation formula together with (6.7), this implies that $(f_\delta^{\alpha_0 m})'(v) > c$ for any $v \in x_m(y_m)'$ and since this is true for any $\alpha_0$ satisfying (6.8), the same estimate holds for $(f_\delta^m)'(v)$.

By Lemma 6.2, $v = -V_m(y_m) + O(\mu_m)$ which implies that $(f_\delta^m)'(-V_m(y)) > c$ as promised. A similar argument shows that $(f_\delta^m)'(V_m) < -c$ which concludes the proof of Lemma 6.3. □

Recall that by (6.5),
$$\{f_\delta = a\} \subset B(o, \delta^3/4) \setminus \bar{B}(o, \nu\delta^3)$$

Therefore
$$\{f_\delta^m = a\} \subset B(x_m, \delta^3/4) \setminus \bar{B}(x_m, \nu\delta^3)$$

If $m$ is sufficiently large. Since $V_m$ is gradient-like for both $f_\delta^m$ and $d(\cdot, x_m)$, a standard flow argument now implies that

$$\{f_\delta^m \geq a\} \text{ is homeomorphic to } \bar{B}(x_m, \delta^3/4) \tag{6.11}$$

for all sufficiently large $m$.

Let $s = n - d$ be the codimension of the collapse. The next lemma is the key ingredient in the proof of Theorem 1.1.

**Lemma 6.4.** *The function $f_\delta^m$ belongs to $\mathcal{T}(s+1)$ in $\bar{B}(x_m, \delta^3/2)$ for any sufficiently large $m$.*

*Proof.* By Remark 3.2, it is enough to show that $f_\delta^{\alpha m} \in \mathcal{T}(s+1)$ for any $\alpha$ and all sufficiently large $m$. Let us fix any $\alpha \in \{1, \ldots, N\}$.

As was explained in section 5, distance functions $d(\cdot, q_m^{\alpha\beta})$ are $-\lambda$-concave in $\bar{B}(x_m, \delta^3/2)$ for some universal positive constant $\lambda$ depending only on the lower curvature bound $k$.

Let $y \in \bar{B}(x_m, \delta^3/2)$.

To prove that $f_\delta^{\alpha m} \in \mathcal{T}(s)$, it is enough to show that for any $\mu_m$-orthonormal frame $v_1, \ldots, v_s \in T_y M_m$, the following estimate holds

$$\sum_{i=1}^{s+1} C f_\delta^{\alpha m}(x, v_i) \leq -\lambda \tag{6.12}$$

Suppose that estimate (6.12) is false and that for some $\mu_m$-orthonormal frame $v_1, \ldots, v_{s+1} \in T_y M_m$ we have

$$\sum_{i=1}^{s+1} C f_\delta^{\alpha m}(x, v_i) > -\lambda \tag{6.13}$$

By Remark 5.2, $f_\delta^{\alpha m}$ is $-\lambda$-concave in $\bar{B}(x_m, \delta/2)$ for *any* choice of $\delta$. Therefore $C f_\delta^m(x, v_i) \leq \lambda$ for any $i = 1, \ldots, s$. Combined with (6.13) this implies that $C f_\delta^{\alpha m}(x, v_i) \geq -(s+1)\lambda$ for every $i = 1, \ldots, s+1$.

Applying Toponogov comparison to the triangle $\triangle y q_m^{\alpha\beta_1} q_m^{\alpha\beta_2}$ we see that $\sphericalangle q_m^{\alpha\beta_1}(y)' q_m^{\alpha\beta_2}(y)' \geq c\delta$ for any $\beta_1 \neq \beta_2$.

Let $\mathcal{A} = \{1, \ldots, N_\alpha\}$. For each $i = 1, \ldots, s+1$, let $\mathcal{A}'_i = \{\beta \in \mathcal{A} | |\cos \sphericalangle \xi v_i| \leq c\delta/4$ for any $\xi \in q_m^{\alpha\beta}(y)'\}$.

Let $\mathcal{A}' = \cap_{i=1}^{s+1} \mathcal{A}'_i$.

Since $v_1, \ldots, v_{s+1}$ are $\mu_m$-orthonormal, this implies that there exists a $d-2$-dimensional totally geodesic sphere $S \subset \Sigma_y M_m$ such that the set $\{q_m^{\alpha\beta}(y)' | \beta \in \mathcal{A}'\}$ lies in the $(c\delta/4 + \mu_m)$- neighborhood of $S$.

For large $m$ we can assume that $(c\delta/4 + \mu_m) \leq c\delta/3$ and since $\sphericalangle q_m^{\alpha\beta_1}(y)'q_m^{\alpha\beta_2}(y)' \geq c\delta$ for any $\beta_1 \neq \beta_2$, a standard volume comparison argument implies that $|\mathcal{A}'| \leq c_3/\delta^{d-2}$.

Since $N \leq c/\delta^{d-1}$ by (6.3), this implies that

$$\frac{|\mathcal{A}'|}{N} \leq c\delta \tag{6.14}$$

We will give separate estimates for $\sum_{i=1}^{s+1} C\phi_d(d(\cdot, q_m^{\alpha\beta}))(x; v_i)$ for $\beta \in \mathcal{A}'$ and for $\beta \notin \mathcal{A}'$.

**Claim 1:** If $\beta \in \mathcal{A}'$ then

$$\sum_{i=1}^{s+1} C\phi_\delta(d(\cdot, q_m^{\alpha\beta}))(y; v_i) \leq (s+1)\lambda \tag{6.15}$$

This follows directly from $-\lambda$-concavity of $\phi_d(d(\cdot, q_m^{\alpha\beta}))$.

**Claim 2:** If $\beta \notin \mathcal{A}'$ then

$$\sum_{i=1}^{s+1} C\phi_\delta(d(\cdot, q_m^{\alpha\beta}))(y; v_i) \leq -\frac{c}{\delta} \tag{6.16}$$

The proof of Claim 2 is essentially the same as the proof of [Kap99, Lemma 4.2] and thus we will skip some of the technical details.

Let us assume for simplicity that $y$ is not a cut point for any of $d(\cdot, q_m^{\alpha\beta})$ so that all the functions involved are actually smooth near $y$.

Let $v \in T_y M_m$ be a unit vector and $\gamma_v(t)$ be a geodesic through $y$ such that $\gamma_v'(0) = v$. Let $f_m^{\alpha\beta}(t) = d(\gamma_v(t), q_m^{\alpha\beta})$.

Then

$$C\phi_\delta(d(\cdot, q_m^{\alpha\beta}))(y; v) = \phi_d''(f_m^{\alpha\beta}(0))((f_m^{\alpha\beta})'(0))^2 + \phi_\delta'(f_m^{\alpha\beta}(0))((f_m^{\alpha\beta})''(0))$$

Observe that $1/2 \leq \phi_\delta'' \leq 1, \phi_\delta'' = -1/2\delta^3$ by the construction of $\phi_\delta$. We also know that $(f_m^{\alpha\beta})'' \leq \lambda$ by the $-\lambda$-concavity of the distance functions and therefore

$$C\phi_\delta(d(\cdot, q_m^{\alpha\beta}))(y; v) \leq c_3 - 1/2\delta^3((f_m^{\alpha\beta})'(0))^2 \tag{6.17}$$

If $\beta \notin \mathcal{A}'$ then by the definition of $\mathcal{A}'$, there exists an $i$ such that $|\cos\sphericalangle q_m^{\alpha\beta}(y)'v_i| \geq c\delta/4$ and therefore, by the first variation formula, $|f_m^{\alpha'}(0)| \geq c\delta/4$ which by (6.17) implies that

$$C\phi_\delta(d(\cdot, q_m^{\alpha\beta}))(y; v_i) \leq c_3 - c_4\delta^3/\delta^2 \leq -c_5/\delta \tag{6.18}$$

Because of $-\lambda$-concavity of the distance functions, for $j \neq i$ we still have that

$$C\phi_\delta(d(\cdot, q_m^{\alpha\beta}))(y; v_j) \leq \lambda \tag{6.19}$$

and therefore

$$\sum_{i=1}^{s+1} C\phi_d(d(\cdot, q_m^{\alpha\beta}))(y; v_i) \le s\lambda - c_5/\delta \le -c_6/\delta \tag{6.20}$$

which concludes the proof of Claim 2 under the extra assumption that $y$ is not a cut point for any of $d(\cdot, q_m^{\alpha\beta})$. The proof of Claim 2 in general is a rather tedious and mostly unilluminating exercise in using the discrete approximation for the formula

$$\phi(h)'' = \phi''(h)(h')^2 + \phi'(h)h''$$

and thus is left to the reader (Also see the proof of [Kap99, Lemma 4.2] where this computation is carried out in detail).

Using claim 1, claim 2 and estimate (6.14) we obtain

$$\sum_{i=1}^{s+1} Cf_\delta^{\alpha m}(x, v_i) \le \frac{|\mathcal{A}'|}{N}((s+1)\lambda) - (1 - \frac{|\mathcal{A}'|}{N})\frac{c_6}{\delta} \le c_7\delta - c_8/\delta \le -c_9/\delta \tag{6.21}$$

Finally, since $c_9$ is independent of $\delta$, we can assume that $\delta$ was chosen to be sufficiently small so that $-c_9/\delta < -\lambda$.

This concludes the proof of Lemma 6.4. □

The first variation formula shows that $f_\delta^m$ satisfies condition (i) of Lemma 4.2. By Lemma 6.4, it also satisfies condition (ii) of Lemma 4.2. Therefore we can apply Lemma 4.2 to $f_\delta^m$ and conclude that $\{f_\delta^m \ge a\}$ satisfies conditions (1)-(3) of Theorem 6.1. Finally, since by (6.11), $\{f_\delta^m \ge a\}$ is homeomorphic to $\bar{B}(x_m, \delta^3/4)$ for all large $m$, we see that the same is true for $\bar{B}(x_m, \delta^3/4)$ as well.

This concludes the proof of Theorem 6.1.

□

**Remark 6.5.** Since the proof of Theorem 6.1 is local on $X$, the theorem remains true for pointed Gromov-Hausdorff convergence.

**Remark 6.6.** When the limit space $X$ in the settings of Theorem 6.1 is a Riemannian manifold, then Yamaguchi's fibration theorem implies that $\bar{B}(x_m, r_0)$ fibers over $\bar{B}(x_0, r_0)$ with the fiber $F_m$ being a closed topological manifold. When $r_0$ is sufficiently small, $\bar{B}(x_0, r_0)$ is contractible and hence, $\bar{B}(x_m, r_0)$ is homotopy equivalent to a closed manifold of expected dimension $n - \dim X$. In other words, if $X$ is smooth then the $CW$ complex provided by Theorem 6.1 can be chosen to be a *closed manifold*.

The author suspects that this remains true for an arbitrary limit space except that the dimension of that manifold can be strictly smaller than $n - \dim X$.

**Remark 6.7.** Observe that for any fixed positive $r < r_0$ provided by Theorem 6.1, functions $d(\cdot, x_m)$ have no critical points in the annuli $\bar{B}(x_m, r_0) \setminus B(x_m, r)$ and therefore the statement of the theorem also holds for $\bar{B}(x_m, r)$ once $m$ is sufficiently large.

**Corollary 6.8.** *Under assumptions of Theorem 6.1, let $\hat{S}(x_m, r_0)$ be the orientation cover of $S(x_m, r_0)$. Then for all sufficiently large $m$, $H^i(\hat{S}(x_m, r_0), A) = 0$ for $n - \dim X < i < \dim X - 1$ and any ring of coefficients $A$. Moreover, if $\dim X \geq 3$, then the same is true for any finite oriented cover of $S(x_m, r_0)$.*

*Proof.* Let $\hat{\bar{B}}(x_m, r_0)$ be the orientation cover of $\bar{B}(x_m, r_0)$. Then $\partial \hat{\bar{B}}(x_m, r_0) = \hat{S}(x_m, r_0)$. By Theorem 6.1 $\bar{B}(x_m, r_0)$ (and hence $\hat{\bar{B}}(x_m, r_0)$) has the homotopy type of a $CW$ complex of $\dim \leq n - \dim X$. Therefore, by Lemma 4.1, $H^i(\hat{S}(x_m, r_0), A) = 0$ for $n - \dim X < i < \dim X - 1$ and any $A$.

Let us suppose that $\dim X \geq 3$. By Theorem 6.1, $\bar{B}(x_m, r_0)$ has the homotopy type of $S(x_m, r_0)$ with a finite number of cells of dimension $\geq \dim X \geq 3$ attached to it. Therefore, the inclusion $S(x_m, r_0) \hookrightarrow \bar{B}(x_m, r_0)$ induces an isomorphism on $\pi_1$. Hence, for any subgroup $\Gamma \subset \pi_1(S(x_m, r_0))$, the corresponding cover $S^\Gamma(x_m, r_0)$ bounds the corresponding cover $\bar{B}^\Gamma(x_m, r_0)$. Since $\bar{B}^\Gamma(x_m, r_0)$ still has the homotopy type of a $CW$ complex of $\dim \leq n - \dim X$, Lemma 4.1 immediately yields the conclusion of the Corollary. $\square$

**Remark 6.9.** The same argument as in the proof of Corollary 6.8 shows that $H^i(S(x_m, r_0), \mathbb{Z}_2) = 0$ for $n - \dim X < i < \dim X - 1$.

# 7 Collapsing to spaces with isolated singularities

The purpose of this section is to prove Theorem 1.5 stated in the introduction. Our main technical statement is the following

**Theorem 7.1.** *Suppose $M_m^n \xrightarrow[m \to \infty]{G-H} X$ where $M_m^n$ is a sequence of $n$-dimensional Riemannian manifolds with $\sec \geq k$ for some $n, k$. Suppose there exists $x_0 \in X$ such that $\Sigma = \Sigma_{x_0} X$ is a closed Riemannian manifold. Then there exists $r_0 = r_0(x_0)$ such that for any $M_m \ni x_m \to x_0$ we have:*

*For any sufficiently large $m$, there exists a topological fiber bundle $F_m \hookrightarrow S(x_m, r_0) \to \Sigma_{x_0} X$ such that*

(1) *$F_m$ and $S(x_m, r_0)$ are closed topological manifolds;*

(2) *$F_m$ is connected;*

(3) *$\pi_1(F)$ is virtually nilpotent;*

(4) *$H^i(\hat{S}(x_m, r_0), A) = 0$ for $\dim F_m < i < n - 1 - \dim F_m$ and any ring of coefficients $A$ where $\hat{S}(x_m, r_0)$ is the orientation cover of $S(x_m, r_0)$; Moreover, if $\dim \Sigma \geq 2$ then the same is true for any oriented finite cover of $S(x_m, r_0)$.*

*Proof.* Let $(M_m^n, x_m)$ be a sequence of manifolds with sectional curvatures bounded below by $k$ converging to $(X, x_0)$.

By the same rescaling argument as in the proof of Theorem 1.1, we can assume that

$$\sec(M_m) \geq -1, X = C\Sigma \text{ and } (M_m, x_m) \xrightarrow[m \to \infty]{G-H} (C\Sigma, o) \qquad (7.22)$$

By Theorem 6.1, there exists an $r_0 > 0$ such that for any sufficiently large $m$, the ball $\bar{B}(x_m, r_0)$ is a thickening of an $s$-dimensional $CW$ complex.

Let $V_m$ be the almost radial vector field on $\bar{B}(x_m, 2r_0) \backslash B(x_m, r_0/2)$ constructed in Lemma 6.2.

Observe that $X$ is a smooth Riemannian manifold of $\sec \geq 0$ away from the vertex $o$. Therefore, by Yamaguchi's fibration theorem, for any sufficiently large $m$ there exists an almost Riemannian submersion with connected fibers $\pi_m \colon U_m \to B(o, 4r_0) \backslash \bar{B}(o, r_0/4)$ where $U_m$ is an open subset of $M_m$ satisfying

$$B(x_m, 2r_0) \backslash B(x_m, r_0/2) \Subset U_m \Subset B(x_m, 8r_0) \backslash B(x_m, r_0/8)$$

Let $V$ be the gradient field of $d(\cdot, o)$ on $X \backslash \{o\}$.

By [Yam91, Lemma 2.8] and Lemma 6.2,

$$\frac{|d\pi_m(V_m(y)) - V(\pi_m(y))|}{|V(\pi_m(y))|} \leq O(\mu_m)$$

for all $y \in U_m$. Therefore $V_m$ is almost perpendicular to the level sets $\pi_m^{-1}(S(o,r))$ for $r_0/2 \leq r \leq 2r_0$. Hence, using the same flow argument as before we obtain that $\pi_m^{-1}(S(o, r_0))$ is homeomorphic to $S(x_m, r_0)$. Thus $S(x_m, r_0)$ fibers over $S(o, r_0)$ which is obviously homeomorphic to $\Sigma$. Let $\pi_m \colon S(x_m, r_0) \to \Sigma$ be the above fibration. To check conditions (2) and (3) observe that by [FY92], the fiber $F$ of $\pi_m$ is connected and has virtually nilpotent fundamental group.

Condition (4) is an immediate consequence of Corollary 6.8.

This concludes the proof of Theorem 7.1. □

**Remark 7.2.** The same argument as in the proof of Theorem 7.1 shows that $H^i(S(x_m, r_0), \mathbb{Z}_2) = 0$ for any $i$ satisfying $\dim F_m < i < n - \dim F_m - 1$.

For a given positively curved manifold $\Sigma$, one can often check that the total space of a bundle $P \to \Sigma$ can never satisfy conditions (1)-(4) of Theorem 7.1 if the dimension of the fiber is too small.

When applied to $\Sigma = F_4/Spin(8), CaP^2$ and $HP^n$ this yields the conclusion of Theorem 1.5:

*Proof of Theorem 1.5 (a).* Suppose $X \in \bar{\mathcal{M}}_k^n$ and $x_0 \in X$ is such that $\Sigma_{x_0} X$ is diffeomorphic to $\Sigma^{24} = F_4/Spin(8)$.

Let $s < 15$ be a positive integer.

**Claim**: For *any* bundle $F \hookrightarrow P^{24+s} \to \Sigma^{24}$ satisfying conditions

1. $\pi_1(F)$ is nilpotent;

2. $F$ is a closed topological manifold;

there exists an $i \in \{dimF+1,\ldots,23\}$ such that $H^i(P,\mathbb{Z}_2) \neq 0$.

Without loss of generality we can assume that $F$ is connected.

Let us look at the $\mathbb{Z}_2$- Serre spectral sequence of the fibration $F \hookrightarrow P^{24+s} \to \Sigma$. It is elementary to check that the nontrivial $\mathbb{Z}_2$ Betti numbers of $\Sigma$ are as follows: $b_0 = b_{24} = 1$, $b_8 = b_{16} = 2$.

If $s < 7$ then the spectral sequence collapses on the $E_2$ term for degree reasons and therefore, $H^8(N, \mathbb{Z}_2) \neq 0$.

Now let $s = 7$. Since $E_7^{0,7} = H^7(F, \mathbb{Z}_2) \cong \mathbb{Z}_2$, and $E_2^{8,0} \cong \mathbb{Z}_2^2$ we see that $\dim_{\mathbb{Z}_2} E_\infty^{8,0} = \dim_{\mathbb{Z}_2} E_8^{8,0} \geq 2 - 1 = 1$. Hence we once again conclude that $H^8(P, \mathbb{Z}_2) \neq 0$.

Next look at the case $8 \leq s \leq 13$. For degree reasons we have that $d_r = 0$ for $1 < r < 7$. Consider $d_7|_{E_7^{0,s}} \colon E_7^{0,s} \cong H^0(\Sigma, \mathbb{Z}_2) \otimes H^s(F, \mathbb{Z}_2) \to E_7^{8,s-7} \cong H^8(\Sigma, \mathbb{Z}_2) \otimes H^{s-7}(F, \mathbb{Z}_2)$.

Since $b_8(\Sigma) = 2$ and $b_s(F) = 1$, this map is either identically zero (if $b_{s-7} = 0$) or not onto (if $b_{s-7} > 0$).

In the former case we conclude by the multiplicativity of the spectral sequence, that $d_r|_{E_r^{8,s}} = 0$ for any $r > 1$. Hence $E_\infty^{8,s} \cong \mathbb{Z}_2^2$ and therefore $H^{s+8}(P, \mathbb{Z}_2) \neq 0$.

If $d_7|_{E_7^{0,s}} \colon E_7^{0,s} \to E_7^{8,s-7}$ is not onto then $E_8^{8,s-7} \cong E_7^{8,s-7}/d_7(E_7^{0,s}) \neq 0$. For degree reasons we have that $E_8^{8,s-7} \cong E_\infty^{8,s-7}$ and therefore, $H^{s+1}(P, \mathbb{Z}_2) \neq 0$.

Let us finally consider the case $s = 14$. Then as in the previous case $d_7|_{E_7^{0,s}} \colon E_7^{0,s} \to E_7^{8,s-7}$ is either zero or not onto. If it is zero then the same argument as before shows that $H^{s+8}(P, \mathbb{Z}_2) \neq 0$.

If $d_7|_{E_7^{0,14}} \colon E_7^{0,14} \to E_7^{8,7}$ is not onto then $E_7^{8,7} \cong H^8(\Sigma, \mathbb{Z}_2) \otimes H^7(F, \mathbb{Z}_2) \neq 0$. Therefore, $b_7(F) \neq 0$. By Poincare duality, the cupproduct on $H^7(F, \mathbb{Z}_2)$ is nondegenerate; therefore $b_7(F)$ is even and hence $\dim_{\mathbb{Z}_2}(E_7^{8,7}) \geq 4$.

Since $\dim_{\mathbb{Z}_2}(E_7^{0,14}) = 1$ and $\dim_{\mathbb{Z}_2}(E_7^{16,0}) = 2$ we see that $\dim_{\mathbb{Z}_2}(E_8^{8,7}) \geq 4 - 2 - 1 = 1$. For degree reasons $E_8^{8,7} = E_\infty^{8,7}$ and therefore $H^{15}(P, \mathbb{Z}_2) \neq 0$. Thus our claim is proved and therefore, by Theorem 7.1, the collapsing codimension $n - 25$ must be least 15.

This concludes the proof of Theorem 1.5 (a). □

*Proof of Theorem 1.5(b).* Suppose $X \in \bar{\mathcal{M}}_k^{17+s}$.

By Theorem 7.1, there exists a bundle $F \hookrightarrow P \to CaP^2$ such that $F$ and $P$ are connected oriented manifolds, $\pi_1(F)$ is virtually nilpotent, $\dim F = s$ and $H^i(\hat{P}, A) = 0$ for any $A$, $s < i < \dim P - s$ and any finite cover $\hat{P} \to P$.

If $s < 7$ then the spectral sequence of this fibration collapses on the $E_2$ term for degree reasons and therefore, $H^8(P) \neq 0$. Thus $s \geq 7$.

Now suppose $s = 7$.

Suppose $\pi_1(F) \neq 0$.

Since $CaP^2$ is 2-connected, the inclusion $F \to P$ is an isomorphism on $\pi_1$. Therefore, after passing to a finite cover of $P$, we can assume that $\pi_1(F)$ is nilpotent.

Hence, there exists an integer $p > 1$ such that $H^1(F, \mathbb{Z}_p) \neq 0$. By Poincare duality, there exist $\alpha \in H^1(F, \mathbb{Z}_p), \beta \in H^6(F, \mathbb{Z}_p)$ such that $\alpha \cup \beta \neq 0$. Let us look at the $\mathbb{Z}_p$ cohomology spectral sequence of $F \hookrightarrow P \to CaP^2$. For degree reasons, $d_j = 0$ for $1 < j < 7$. By multiplicutivity of the spectral sequence, $d_7(\alpha \cup \beta) = 0$ and hence $d_7|_{E_7^{0,7}} = 0$. Therefore $H^8(P, \mathbb{Z}_p) \neq 0$.

Now suppose that $\pi_1(F) = 0$.

We consider two possibilities:

**Case 1.** There exists an $i$ satisfying $1 < i < 6$ such that $H^i(P) \neq 0$. Choose the minimal $i$ satisfying this condition. Then there exists a $p$ such that $H^i(P, \mathbb{Z}_p) \neq 0$. The same Poincare duality argument as above now implies that $H^8(P, \mathbb{Z}_p) \neq 0$.

**Case 2.** $H^i(F) = 0$ for any $0 < i < 7$. Since $F$ is simply connected, by Poincare conjecture, this implies that $F$ is homeomorphic to $S^7$. Let us look at the $\mathbb{Z}$-cohomology spectral sequence of $F \hookrightarrow P \to CaP^2$. If $d_7|_{E_7^{0,7}} \colon E_7^{0,7} \cong \mathbb{Z} \to E_7^{8,0} \cong \mathbb{Z}$ is not an isomorphism then $H^8(P) \neq 0$. If $d_7|_{E_7^{0,7}} \colon E_7^{0,7} \to E_7^{8,0}$ is an isomorphism then $H^i(P) = 0$ for any $0 < i < 23$. Therefore $P$ is homeomorphic to $S^{23}$ by the Poincare conjecture.

However, according to [Bro63], there is no bundle with the total space $S^{23}$ and the fiber homeomorphic to $S^7$. Thus this case is also impossible and the codimension of collapse $s$ can not be equal to 7. Therefore $s \geq 8$ as claimed. $\square$

*Proof of Theorem 1.5(c).* It is easy to see that if $F \hookrightarrow P \to \mathbb{HP}^n$ is a bundle with $\dim F \leq 2$, then the Serre spectral sequence of this fibration collapses on the $E_2$ term and therefore $H^4(P) \neq 0$. By Theorem 7.1, this implies that the minimal collapsing codimension of $X$ is at least 3.

On the other hand, the same observation of Yamaguchi which we used to compute the minimal collapsing codimension of $\mathbb{SCP}^n$, shows that there exists a sequence of metrics with $\sec \geq 0$ on $S(S^{4n+3}) = S^{4n+4}$ converging to $\mathbb{SHP}^n$. Therefore, the minimal collapsing codimension of $\mathbb{SHP}^n$ is equal to 3. $\square$

*Proof of Theorem 1.7.* Suppose $M_m^{4n+4} \xrightarrow[m \to \infty]{G-H} X = S\mathbb{HP}^n$ where $n \geq 2$.

Let $x_0$ and $y_0$ be the north and the south poles of $X$. Let $r_{x_0}, r_{y_0}$ be the radii provided by Theorem 6.1. Let $M_m \ni x_m \to x_0$, $M_m \ni y_m \to y_0$. Since $d(\cdot, x_0)$ and $d(\cdot, y_0)$ have no critical points in $X \backslash (B(x_0, r_{x_0}) \cup B(y_0, r_{y_0}))$, a standard critical point argument shows that $d(\cdot, x_m)$ and $d(\cdot, y_m)$ have no critical points in $M_m \backslash (B(x_m, r_{x_0}) \cup B(y_m, r_{x_0}))$, and moreover,

$$M_m \backslash (B(x_m, r_{x_0}) \cup B(y_m, r_{x_0})) \text{ is homeomorphic to } S(x_m, r_{x_0}) \times [0, 1] \qquad (7.23)$$

By the same argument as in the proof of part (4) of Theorem 7.1, we see that inclusions $S(x_m, r_{x_0}) \hookrightarrow \bar{B}(x_m, r_{x_0})$, $S(y_m, r_{y_0}) \hookrightarrow \bar{B}(y_m, r_{y_0})$ are isomorphisms on $\pi_1$ and therefore the inclusion $S(x_m, r_{x_0}) \hookrightarrow M_m$ is also an isomorphism on $\pi_1$ if $m$ is sufficiently large.

We will first show that the universal cover of $M_m$ is homeomorphic to a sphere.

By Theorem 7.1, $S(x_m, r_{x_0})$ fibers over $\mathbb{HP}^n$ with the fiber $F_m$ being a closed 3-manifold with virtually nilpotent fundamental group. Since $\mathbb{HP}^n$ is 2-connected, the inclusion $F_m \hookrightarrow S(x_m, r_{x_0})$ is an isomorphism on $\pi_1$ and by above, the same is true for $F_m \hookrightarrow M_m$.

Let $\Gamma \leq \pi_1$ be a nilpotent subgroup of $\pi_1(F_m)$ of finite index such that the corresponding cover $F^\Gamma$ of $F_m$ is orientable.

Let $M_m^\Gamma, \bar{B}^\Gamma(x_m, r_{x_0}), S^\Gamma(x_m, r_{x_0}), \bar{B}^\Gamma(y_m, r_{y_0})$ be the corresponding covers of $M_m$, $\bar{B}(x_m, r_{x_0})$, $S(x_m, r_{x_0})$ and $\bar{B}(y_m, r_{y_0})$ respectively.

By above,
$$M_m^\Gamma \text{ is homeomorphic to } \bar{B}^\Gamma(x_m, r_{x_0}) \cup_{f^\Gamma} \bar{B}^\Gamma(y_m, r_{y_0}) \tag{7.24}$$
where $f^\Gamma$ is the natural identification of $S^\Gamma(x_m, r_{x_0})$ and $S^\Gamma(y_m, r_{y_0})$ induced by (7.23).

Suppose that $\Gamma \neq 1$.

Since $\Gamma$ is nilpotent, $H^1(F, \mathbb{Z}_p) \neq 0$ for some $p > 1$. By Poincare duality, there exist $a \in H^1(F, \mathbb{Z}_p), b \in H^2(F, \mathbb{Z}_p)$ such that $a \cup b \neq 0$. By looking at the $\mathbb{Z}_p$ cohomology spectral sequence of $F_m^\Gamma \to S^\Gamma(x_m, r_{x_0}) \to \mathbb{HP}^n$ we see that it collapses on the $E_2$ term and thus $H^4(S^\Gamma(x_m, r_{x_0}), \mathbb{Z}_p) \neq 0$. This is impossible by Theorem 7.1 and hence $\Gamma = 1$.

Therefore all the spaces $M_m^\Gamma, F_m^\Gamma \bar{B}^\Gamma(x_m, r_{x_0}), S^\Gamma(x_m, r_{x_0}), \bar{B}^\Gamma(y_m, r_{y_0})$ are simply connected. In particular, $F_m^\Gamma$ is a homotopy 3-sphere.

By looking at the $\mathbb{Z}$ cohomology spectral sequence of $F_m^\Gamma \to S^\Gamma(x_m, r_{x_0}) \to \mathbb{HP}^n$ we see that if the Euler class of this fibration is not a generator of $H^4(\mathbb{HP}^n)$, then $H^4(S^\Gamma(x_m, r_{x_0})) \neq 0$. Since we know that this is impossible, the Euler class is a generator of $H^4(\mathbb{HP}^n)$ and therefore, $S^\Gamma(x_m, r_{x_0})$ is an integral homology sphere. By above, $S^\Gamma(x_m, r_{x_0})$ is simply connected and hence it is homeomorphic to $S^{4n+3}$ by the Poincare conjecture.

This means that $\bar{B}^\Gamma(x_m, r_{x_0})$ is contractible. Indeed, we already know that it is simply connected and by Theorem 1.1, it has the homotopy type of a 3-dimensional complex. Thus $H^i(\bar{B}^\Gamma(x_m, r_{x_0})) = H_i(\bar{B}^\Gamma(x_m, r_{x_0})) = 0$ for any $i > 3$.

Look at the long exact cohomology sequence of the pair $(\bar{B}^\Gamma(x_m, r_{x_0}), S^\Gamma(x_m, r_{x_0}))$:

$$\to H^i(\bar{B}^\Gamma(x_m, r_{x_0}), S^\Gamma(x_m, r_{x_0})) \to H^i(\bar{B}^\Gamma(x_m, r_{x_0})) \to H^i(S^\Gamma(x_m, r_{x_0})) \to$$

When $i = 2$ or $3$ we see that $H^i(S^\Gamma(x_m, r_{x_0})) = H^i(S^{4n+3}) = 0$. Also, by Poincare duality, $H^i(\bar{B}^\Gamma(x_m, r_{x_0}), S^\Gamma(x_m, r_{x_0})) = H_{4n+3-i}(\bar{B}^\Gamma(x_m, r_{x_0})) = 0$. By the long exact sequence above, this immediately yields that $H^2(\bar{B}^\Gamma(x_m, r_{x_0})) = H^3(\bar{B}^\Gamma(x_m, r_{x_0})) = 0$. Thus $\bar{B}^\Gamma(x_m, r_{x_0})$ is simply connected and $H^i(\bar{B}^\Gamma(x_m, r_{x_0})) = 0$ for any $i > 0$. By Hurewitz Theorem, this implies that $\bar{B}^\Gamma(x_m, r_{x_0})$ is contractible. We also know that its boundary is homeomorphic to a sphere. The same is true for $\bar{B}^\Gamma(y_m, r_{y_0})$ which by (7.24) implies that $M_m^\Gamma$ is a simply connected homology sphere. By Poincare conjecture, $M_m^\Gamma$ must be homeomorphic to $S^{4n+4}$.

It is well-known [Bre72] that the only nontrivial group that can act freely on an even dimensional sphere is $Z_2$. Thus $\pi_1(M_m)$ is either trivial or is isomorphic to $Z_2$. We claim that the latter case is impossible. Indeed, if $\pi_1(M_m) \cong \mathbb{Z}_2$, then $\pi_1(F_m) \cong \pi_1(M_m)$ is also isomorphic to $\mathbb{Z}_2$.

Looking at the $\mathbb{Z}_2$ cohomology spectral sequence of the bundle $F_m \to S(x_m, r_{x_0}) \to \mathbb{HP}^n$, the same argument as before shows that $H^4(S(x_m, r_{x_0}), \mathbb{Z}_2) \neq 0$ which is impossible by Theorem 7.1. Thus, $\pi_1(M_m) = 1$ and hence, $M_m = M_m^\Gamma$ is homeomorphic to $S^{4n+4}$ as claimed.

$\square$

**Remark 7.3.** It would be interesting to see whether the conclusion of Theorem 1.7 can be improved to show that a manifold $M^{4n+4}$ with $n > 1, \sec(M) \geq k$ sufficiently close to $S\mathbb{HP}^n$ must be *diffeomorphic* to $S^{4n+4}$.

It would also be interesting to see if Theorem 1.7 remains true if the assumption $\sec(M) \geq k$ is replaced by the weaker one $Ricc(M) \geq (n-1)k$. The author believes that this is most likely false.

**Remark 7.4.** Same proof as in Theorem 1.7 shows that if $M_m^{2n+2} \xrightarrow[m \to \infty]{G-H} S\mathbb{CP}^n$ where $\sec(M_m) \geq k$ and $n > 1$ then $M_m^{2n+2}$ is homeomorphic to $S^{2n+2}$ for all large $m$.

## 8 Concluding remarks

Using Theorem 7.1 one can obtain nontrivial bounds for the minimal collapsing codimensions of spherical suspensions or cones over other positively curved manifolds such as Eschenburg and Bazaikin manifolds.

Unfortunately, due to the lack of examples of positively curved manifolds, Theorem 7.1 does not produce examples of Alexandrov spaces with arbitrary large minimal collapsing codimensions. However, the author suspects that it might be possible to show that for any positively curved $\Sigma$ different from a sphere, the minimal collapsing codimension of $\underbrace{\Sigma * \Sigma * \ldots * \Sigma}_{l}$ grows at least linearly in $l$.

Using the same ideas as in the proof of Theorem 1.5, it might also be possible to show that Theorem 1.5 remains true for $X \times M$ where $M$ is any Riemannian manifold (i.e that if, say, $X$ has a point $x$ with $\Sigma_x X$ diffeomorphic to $F_4/Spin(8)$, then $X \times M \notin \bar{\mathcal{M}}_k^{\dim(X \times M)+s}$ for any $k \in \mathbb{R}$ and $s < 15$ ).